\documentclass[12pt]{article}
\usepackage{amssymb,amsfonts,amsmath, psfrag, cite}
\usepackage{graphicx,colordvi}

\newtheorem{theo}{Theorem}

\makeatletter \@addtoreset{equation}{section} \@addtoreset{theo}{section}
\makeatother

\def\qed{\hfill \rule{4pt}{7pt}}
\def\pf{\noindent {\it Proof.} }
\def\D{{\Delta}}

\begin{document}

\begin{center}
{\Large A Telescoping Method for Double Summations}

 William Y.C. Chen$^1$,
Qing-Hu Hou$^2$ and Yan-Ping Mu$^3$\\[5pt]
Center for Combinatorics, LPMC \\
Nankai University, Tianjin 300071, P. R. China \\
E-mail: $^1$chen@nankai.edu.cn, $^2$hou@nankai.edu.cn, $^3$myphb@eyou.com

Dedicated to James D. Louck on the Occasion of His Seventy-Fifth
Birthday

\end{center}

\begin{abstract}

We present a method to prove hypergeometric double summation
identities. Given a hypergeometric term $F(n,i,j)$, we aim to find
a difference operator $ L=a_0(n) N^0 + a_1(n) N^1 +\cdots+a_r(n)
N^r $ and rational functions $R_1(n,i,j),R_2(n,i,j)$ such that $ L
F = \Delta_i (R_1 F) + \Delta_j (R_2 F)$. Based on simple
divisibility considerations, we show that the denominators of
$R_1$ and $R_2$ must possess certain factors which can be computed
from $F(n, i,j)$. Using these factors as estimates, we may find
the numerators of $R_1$ and $R_2$ by guessing the upper bounds of
the degrees and solving systems of linear equations. Our method is
valid for the Andrews-Paule identity, Carlitz's identities, the
Ap\'ery-Schmidt-Strehl identity, the Graham-Knuth-Patashnik
identity, and the Petkov\v{s}ek-Wilf-Zeilberger identity.
\end{abstract}



\section{\normalsize Introduction}

This paper is concerned with double summations of hypergeometric terms $F(n, i,
j)$.  A function $F(n,k_1,\ldots,k_m)$ is called a {\it hypergeometric term} if
the quotients
\[
{F(n+1,k_1,\ldots,k_m) \over F(n,k_1,\ldots,k_m)}, \quad {F(n,k_1+1,\ldots,k_m)
\over F(n,k_1,\ldots,k_m)}, \quad \ldots, \quad {F(n,k_1,\ldots,k_m+1) \over
F(n,k_1,\ldots,k_m)}
\]
are rational functions of $n,k_1,\ldots,k_m$. Throughout the
paper, we use $N$ to denote the shift operator with respect to the
variable $n$, given  by $N F(n)=F(n+1)$ and use $\Delta_{x}$ to
denote the difference operator with respect to the variable $x$,
given by $\Delta_x F = F(x+1)-F(x)$. For polynomials $a$ and $b$,
we denote by $\gcd(a,b)$ their monic greatest common divisor. When
we express a rational function as a quotient $p/q$, we always
assume that  $p$ and $q$ are relatively prime unless it is
explicitly stated otherwise.

Zeilberger's algorithm \cite{Gra,Pet-W-Z,Zeil}, also known as the
method of {\it creative telescoping}, is devised for proving
hypergeometric identities of the form
\begin{equation}
\label{F-f} \sum_{k} F(n,k) = f(n),
\end{equation}
where $F(n,k)$ is a hypergeometric term and $f(n)$ is a given
function. This algorithm has been used to deal with multiple sums
by Wilf and Zeilberger \cite{W-Z92}. Given a hypergeometric term
$F(n,k_1,\ldots,k_m)$, the approach of Wilf and Zeilberger tries
to find a linear difference operator $L$ with coefficients being
polynomials in $n$
$$
L=a_0(n) N^0 + a_1(n) N^1 +\cdots+a_r(n) N^r
$$
and rational functions $R_1,\ldots,R_m$ of $n,k_1,\ldots,k_m$ such that
\begin{equation}
\label{L-R} L F=\sum_{l=1}^m \Delta_{k_l} (R_l F).
\end{equation}
As noted by K.~Wegschaider \cite{Weg}, when the boundary
conditions are admissible, Equation \eqref{L-R} leads to a
homogenous recursion for the multi-summations:
\[
L \sum_{k_1,\ldots,k_m} F(n,k_1,\ldots,k_m) = 0.
\]

When $m=1$, $L$ and $R_1$ can be solved by Gosper's algorithm
\cite{Gosp,Pet-W-Z}. S.A.~Abramov, K.O.~Geddes and H.Q.~Le also
provided a lower bound for the order $r$ \cite{AbGL02,AbL02} and
found a faster algorithm \cite{AbL04} compared with Zeilberger's
algorithm. For a survey on recent developments, see \cite{AbCG04}.
For $m \ge 2$,  constructing the denominators of $R_1,\ldots,R_m$
for the Wilf-Zeilberger approach remains an open problem. In a
recent paper \cite{Moh}, M.~Mohammed and D.~Zeilberger used the
denominator of $LF/F$ as estimates of the denominators of $R_i$.
In an alternative approach, Wegschaider generalized Sister
Celine's technique \cite{Weg} to multiple summations, and proved
many double summation identities. A different approach has been
proposed by F.~Chyzak \cite{Chy-Sal,Chy} by finding recursions of
the summation iteratively starting from the inner sum.
C.~Schneider \cite{Sch} presented the Chyzak method  from the
point of view of Karr's difference field theory.

To give a sketch of our approach, we first consider Gosper's
algorithm for bivariate hypergeometric terms. Suppose that
$F(i,j)$ is a hypergeometric term and $p_1/q_1, p_2/q_2$ are
rational functions such that
\[
F(i,j) = \D_i \left( {p_1(i,j) \over q_1(i,j)} F(i,j) \right) + \D_j \left(
{p_2(i,j) \over q_2(i,j)} F(i,j)\right).
\]
We show that under certain hypotheses (Section $2$, $(H1)$--$(H3)$), the
denominators $q_1,q_2$ can be written in the form
\begin{equation}
\label{factor}
\begin{array}{l}
q_1(i,j) = v_1(i)\, v_2(j)\, v_3(i+j)\, v_4(i,j)\, u_1(j)\, u_2(i,j), \\[7pt]
q_2(i,j) = v_1(i)\, v_2(j)\, v_3(i+j)\, v_4(i,j)\, w_1(i)\,
w_2(i,j),
\end{array}
\end{equation}
such that $v_1,v_2,v_4$ and $u_2,w_2$ are bounded in the sense
that they are factors of certain polynomials which can be computed
for a given $F(i,j)$, see Theorem~\ref{Gosper}. Then we apply
these estimates to the telescoping algorithm for double
summations. Suppose that
\[
L F(n,i,j) = \D_i \big( R_1(n,i,j) F(n,i,j) \big) + \D_j \big( R_2(n,i,j)
F(n,i,j) \big),
\]
where
\[
R_1(n,i,j)={1 \over d(n,i,j)} \cdot {f_1(n,i,j) \over g_1(n,i,j)}, \qquad
R_2(n,i,j)={1 \over d(n,i,j)} \cdot {f_2(n,i,j) \over g_2(n,i,j)}
\]
and $d(n,i,j)$ is the denominator of $LF(n,i,j)/F(n,i,j)$. We may
deduce that $g_1,g_2$ can be factored in the form of
\eqref{factor} such that $v_1,v_2,v_4$ and $u_2,w_2$ are bounded,
see Theorem~\ref{Zeil}. Although we do not have the universal
denominators, these bounds can be used to give estimates of the
denominators $g_1$ and $g_2$. Then by further guessing the bounds
of the degrees of the numerators of $R_1$ and $R_2$, we get the
desired difference operator if we are lucky.

Indeed, our approach works quite efficiently for many identities such as the
Andrews-Paule identity, Carlitz's identities, the Ap\'ery-Schmidt-Strehl
identity, the Graham-Knuth-Patashnik identity, and the
Petkov\v{s}ek-Wilf-Zeilberger identity.

\section{\normalsize Denominators in Bivariate Gosper's Algorithm}
For a given bivariate hypergeometric term $F(i,j)$, we give estimates of the
denominators of the rational functions $R_1(i,j), R_2(i,j)$ satisfying
\begin{equation}
\label{hyper-eqn} F(i,j) = \D_i \big( R_1(i,j) F(i,j)\big) + \D_j \big(
R_2(i,j) F(i,j)\big).
\end{equation}

Let
\begin{equation}
\label{reduce}
\begin{array}{ll}
\displaystyle R_1(i,j) = {f_1(i,j) \over g_1(i,j)},
& \displaystyle  R_2(i,j) = {f_2(i,j) \over g_2(i,j)}, \\[15pt]
\displaystyle  {F(i+1,j) \over F(i,j)} = {r_1(i,j) \over s_1(i,j)}, &
\displaystyle {F(i,j+1) \over F(i,j)} = {r_2(i,j) \over s_2(i,j)}.
\end{array}
\end{equation}
Dividing $F(i,j)$ on both sides of \eqref{hyper-eqn} and substituting
\eqref{reduce} into it, we derive that
\begin{equation}
\label{main-eqn} 1 = {r_1(i,j) \over s_1(i,j)} {f_1(i+1,j) \over g_1(i+1,j)}
- {f_1(i,j) \over g_1(i,j)} + {r_2(i,j) \over s_2(i,j)} {f_2(i,j+1) \over g_2(i,j+1)}
- {f_2(i,j) \over g_2(i,j)}.
\end{equation}

Let
\[
u(i,j)=\gcd(s_1(i,j),s_2(i,j)), \qquad
v(i,j)=\gcd(g_1(i,j),g_2(i,j)),
\]
and
\begin{equation}
\label{u-v}
\begin{array}{l}
s_1'(i,j)=s_1(i,j)/u(i,j), \qquad  s_2'(i,j)=s_2(i,j)/u(i,j), \\[5pt]
g_1'(i,j)=g_1(i,j)/v(i,j), \qquad  g_2'(i,j)=g_2(i,j)/v(i,j).
\end{array}
\end{equation}

We find that in many cases we can restrict our attention to  those
$R_1,R_2$ whose denominators $g_1, g_2$ satisfy the following
three hypotheses. We see  that in the proof of the following
theorem, these hypotheses enable us to cancel out unknown factors
from the multiples of $g_1$ and $g_2$ so that we can obtain an
upper bound of $g_1$ and $g_2$. Thus, these hypotheses come
naturally from the requirement of simple divisibility properties.
Moreover, it turns out that these divisibility requirements are
sufficient in many cases to give good estimates for the
denominators $g_1$ and $g_2$. The three hypotheses are as follows:
\begin{itemize}
\item[$(H1)$] Suppose $p(i,j)$ and $p(i+h_1,j+h_2)$ are both
irreducible factors of $g_1(i,j)$ ($g_2(i,j)$, respectively) for
some $h_1,h_2 \in \{-1,0,1\}$. Then they must be coincide.
\item[$(H2)$] $\gcd(g_1'(i,j),v(i,j))=\gcd(g_2'(i,j),v(i,j))=1$.
\item[$(H3)$] For any integers $h_1,h_2 \in \{-1,0,1\}$,
\[
\gcd(g_1'(i+h_1,j+h_2),g_2'(i,j))=1.
\]
\end{itemize}
For example, the following functions satisfy the above hypotheses:
\[
g_1(i,j)=(2n-2i+1)(n-i+1)(j+1)^2, \quad
g_2(i,j)=(2n-2i+1)(n-i+1)(i+1)^2.
\]

 \noindent {\bf Remarks.}
\begin{itemize}
\item[1.] Hypothesis $(H1)$ looks like requiring that $g_1$ and
$g_2$ are shift-free (see Abramov and Petkov\v{s}ek
\cite{Abra02}). However, only the shifts of $\pm 1$ are considered
and shift invariant factors are admissible. For example, we allow
that $g_1(i,j)=(i+1)(i+3)$ or $g_1(i,j)=i+j$.

\item[2.] According to \cite{Ab-Pet},
$\gcd(g_1(i,j),g_1(i+h_1,j+h_2))$ and
$\gcd(g_2(i,j),g_2(i+h_1,j+h_2))$ can factor into integer-linear
factors for $h_1,h_2$ being not both zero.

\item[3.] Hypothesis $(H3)$ is to require that $g'_1/g'_2$ are
shift-reduced (see also \cite{Abra02}) respect to the shifts of
$\pm 1$.
\end{itemize}

Under the above hypotheses, we have
\begin{theo}
\label{Gosper}
The denominators $g_1(i,j),g_2(i,j)$ can be factored into polynomials:
\[
\begin{array}{l}
g_1(i,j) = v_1(i) v_2(j) v_3(i+j) v_4(i,j) u_1(j) u_2(i,j), \\[5pt]
g_2(i,j) = v_1(i) v_2(j) v_3(i+j) v_4(i,j) w_1(i) w_2(i,j),
\end{array}
\]
such that
\begin{align}
& v_1(i) \,|\, r_1(i-1,j) s_2'(i-1,j), \label{v1}\\[5pt]
& v_2(j) \,|\, r_2(i,j-1) s_1'(i,j-1), \label{v2}\\[5pt]
& v_4(i,j) \,|\, \gcd \big( r_1(i-1,j) s_2'(i-1,j), r_2(i,j-1) s_1'(i,j-1) \big), \label{v4}\\[5pt]
& u_2(i,j) \,|\, \gcd \big( s_1(i,j) s_2'(i,j), r_1(i-1,j) s_2'(i-1,j) \big), \label{u2}\\[5pt]
& w_2(i,j) \,|\, \gcd \big( s_2(i,j) s_1'(i,j), r_2(i,j-1)
s_1'(i,j-1) \big). \label{w2} \end{align}
\end{theo}

\pf Substituting \eqref{u-v} into \eqref{main-eqn}, we get
\begin{eqnarray*}
1 & = & {r_1(i,j) \over s'_1(i,j)u(i,j)} {f_1(i+1,j) \over g'_1(i+1,j)v(i+1,j)}
- {f_1(i,j) \over g'_1(i,j)v(i,j)} \\
&& + {r_2(i,j) \over s'_2(i,j)u(i,j)}
{f_2(i,j+1) \over g'_2(i,j+1)v(i,j+1)}
- {f_2(i,j) \over g'_2(i,j)v(i,j)}.
\end{eqnarray*}
That is,
\begin{eqnarray*}
\lefteqn{ s_1(i,j) s_2'(i,j) g_1(i,j) g'_2(i,j) g_1(i+1,j) g_2(i,j+1) } \\
&=& f_1(i+1,j) r_1(i,j) s_2'(i,j) g_1(i,j) g'_2(i,j) g_2(i,j+1) \\
&& - f_1(i,j) s_1(i,j) s_2'(i,j) g'_2(i,j) g_1(i+1,j) g_2(i,j+1) \\
&& + f_2(i,j+1) r_2(i,j) s'_1(i,j) g_1(i,j) g'_2(i,j) g_1(i+1,j) \\
&& - f_2(i,j) s_1(i,j) s_2'(i,j) g'_1(i,j) g_1(i+1,j) g_2(i,j+1).
\end{eqnarray*}

\begin{enumerate}
\item Suppose that $p(i,j)$ is an irreducible factor of $v(i,j)$,
and for some non-negative integer $l$, $p^l \,|\, v$. Note that
$p(i+h_1,j+h_2)$ is also irreducible. Since
\[
\gcd(p(i+1,j),f_1(i+1,j))= \gcd(p(i,j+1),f_2(i,j+1))=1,
\]
we have
\[
p^l(i+1,j)  \,|\,   r_1(i,j) s_2'(i,j) g_1(i,j) g'_2(i,j) g_2(i,j+1)
\]
and
\[
p^l(i,j+1) \,|\,  r_2(i,j) s'_1(i,j) g_1(i,j) g'_2(i,j) g_1(i+1,j).
\]

There are three cases:
\begin{itemize}
\item $p(i,j)$ is a polynomial depending only on $i$. Then $\gcd(p(i+1,j),
g_1(i,j))=1$. Otherwise, by hypothesis $(H1)$ we have that $p(i+1,j)=p(i,j)$ is
independent of $i$, which is a contradiction. Similarly, $\gcd(p(i+1,j),
g_2(i,j))=1$. Since $p(i,j)$ is a polynomial depending only on $i$, we have
\[
\gcd(p(i+1,j), g_2(i,j+1)) = \gcd(p(i+1,j+1),g_2(i,j+1)) = 1.
\]
Hence,
$$p^l(i+1,j) \,|\, r_1(i,j) s_2'(i,j). $$
Let $v_1(i)$ denote the product of all irreducible factors of
$v(i,j)$ that depend only on $i$. Then we have \eqref{v1}.

\item $p(i,j)$ is a polynomial depending only on $j$. The same discussion leads to
$$p^l(i,j+1) \,|\, r_2(i,j) s'_1(i,j). $$
Let $v_2(j)$ denote the product of all irreducible factors of
$v(i,j)$ that depend only on $j$. Then we have \eqref{v2}.

\item $p(i,j)$ is a polynomial depending both on $i$ and on $j$. Then either
\begin{equation}
\label{v3f} p(i+1,j)=p(i,j+1)
\end{equation}
or
\begin{equation}
\label{v4f} \gcd(p(i+1,j),p(i,j+1))=1.
\end{equation}
In the former case, $p(i,j)$ is a polynomial of $i+j$ (see
\cite[Lemma 3]{Ab-Pet} or \cite[Lemma 3.3]{Hou}). For this case we
do not have a bound. We denote by $v_3(i+j)$ the product of all
irreducible factors $p(i,j)$ of $v(i,j)$ that satisfy \eqref{v3f}.
In the later case, by hypothesis $(H1)$, we have
\[
\gcd(p(i+1,j),  g_1(i,j) g'_2(i,j) g_2(i,j+1)) = 1
\]
and
\[
\gcd(p(i,j+1), g_1(i,j) g'_2(i,j) g_1(i+1,j)) = 1.
\]
Thus,
\[
p^l(i,j) \,|\, \gcd \big( r_1(i-1,j) s_2'(i-1,j), r_2(i,j-1) s_1'(i,j-1) \big).
\]
Let $v_4(i,j)$ denote the product of all irreducible factors
$p(i,j)$ of $v(i,j)$ that satisfy \eqref{v4f}. Then we have
\eqref{v4}.
\end{itemize}

\item Suppose $p$ is an irreducible factor of $g_1'$ and $p^l |
g_1'$ for some non-negative integer $l$. If $p(i,j)\,|\,v(i,j+1)$,
then $p(i,j-1)\,|\,v(i,j)$. By hypothesis $(H1)$,
$p(i,j-1)=p(i,j)$, which implies $p(i,j)\,|\,v(i,j)$,
contradicting to hypothesis $(H2)$. Noting further that by
hypothesis $(H3)$, for any $h_1,h_2 \in \{-1,0,1\}$,
\[
\gcd(f_1(i,j),g_1(i,j))=\gcd(g_1'(i,j),g_2'(i+h_1,j+h_2))=1,
\]
we have that
\[
p^l(i,j) \,|\, s_1(i,j) s_2'(i,j) g_1(i+1,j).
\]
If $p(i+1,j)\,|\,v(i,j+1)$, then by hypothesis $(H1)$,
$p(i+1,j-1)=p(i,j)$, which implies $p(i,j)\,|\,v(i,j)$,
contradicting to hypothesis $(H2)$. Therefore, by hypothesis
$(H3)$,
\[
p^l(i+1,j) \,|\, r_1(i,j) s_2'(i,j) g_1(i,j).
\]
There are two cases:
\begin{itemize}
\item $p(i,j)=p(i+1,j)$. Then $p(i,j)$ is a polynomial depending
only on $j$. For this case we also do not have a bound. We denote
by $u_1(j)$ the product of all irreducible factors of $g'_1(i,j)$
that depend only on $j$.

\item $\gcd(p(i,j),p(i+1,j))=1$. Then by hypothesis $(H1)$,
\[
\gcd(p(i,j), g_1(i+1,j)) = \gcd(p(i+1,j),g_1(i,j)) = 1,
\]
and hence,
\[
p^l(i,j) \,|\, \gcd \big( s_1(i,j) s_2'(i,j),  r_1(i-1,j)
s_2'(i-1,j) \big).
\]
Let $u_2(i,j)$ denote the product of all irreducible factors
$p(i,j)$ of $g'_1(i,j)$ such that $\gcd(p(i,j),p(i+1,j))=1$. Then
 we have \eqref{u2}.
\end{itemize}

\item Similarly, suppose $p$ is an irreducible factor of $g_2'$ and $p^l |
g_2'$ for some non-negative integer $l$.  Then either $p(i,j)$ is a polynomial
depending only on $i$ or
\[
p^l(i,j) \,|\, \gcd \big( s_2(i,j) s_1'(i,j),  r_2(i,j-1)
s_1'(i,j-1) \big).
\]
Let $w_1(i)$ denote product of irreducible factors of $g'_2(i,j)$
that depend only on $i$ and $w_2(i,j)$ denote the product of the
rest irreducible factors of $g'_2(i,j)$. Then we have \eqref{w2}.
\qed
\end{enumerate}
 Note that $u_2(i,j)$ have no factors which are
free of $i$ and $w_2(i,j)$ have no factors which are free of $j$.
We will need this property later for the algorithm EstDen.

\section{\normalsize Denominators in Our Telescoping Method}

We are now ready to estimate the denominators of $R_1$ and $R_2$
in our telescoping method.

As in the case of single summations, the telescoping algorithm for double
summations tries to find an operator
\[
L = a_0(n) + a_1(n)N + \cdots a_r(n)N^r
\]
and rational functions $R_1(n,i,j), R_2(n,i,j)$ such that
\begin{equation}
\label{WZ-triple} L F(n,i,j) = \D_i (R_1(n,i,j) F(n,i,j)) + \D_j (R_2(n,i,j)
F(n,i,j)).
\end{equation}

Let
\begin{equation}
\label{IJ-rat} {F(n,i+1,j) \over F(n,i,j)} = {r_1(n,i,j) \over s_1(n,i,j)},
\qquad {F(n,i,j+1) \over F(n,i,j)} = {r_2(n,i,j) \over s_2(n,i,j)},
\end{equation}
and $d(n,i,j)$ be the common denominator of
\[
{F(n+1,i,j) \over F(n,i,j)}, \quad \ldots, \quad {F(n+r,i,j) \over F(n,i,j)}.
\]
Then there exists a polynomial $c(n,i,j)$, not necessarily being coprime to
$d$, such that
\begin{equation}
\label{f/F} {LF(n,i,j) \over F(n,i,j)} = \sum_{l=0}^r a_l(n) {F(n+l,i,j) \over
F(n,i,j)} = {c(n,i,j) \over d(n,i,j)}.
\end{equation}
Note that $c$ is related to the polynomials $a_0,a_1,\ldots,a_r$ but $d$ is
independent of them.

Now, \eqref{WZ-triple} can be written in the form of \eqref{hyper-eqn}:
\[
L F(n,i,j) = \D_i (R_1'(n,i,j) LF(n,i,j)) + \D_j (R_2'(n,i,j) LF(n,i,j)),
\]
where
\[
R'_1(n,i,j) = R_1(n,i,j) {d(n,i,j) \over c(n,i,j)} \quad \mbox{and} \quad
R'_2(n,i,j) = R_2(n,i,j) {d(n,i,j) \over c(n,i,j)}.
\]
This suggests us to assume
\begin{equation}
\label{f-g} R_1(n,i,j)={1 \over d(n,i,j)} {f_1(n,i,j) \over g_1(n,i,j)} \quad
\mbox{and} \quad R_2(n,i,j)={1 \over d(n,i,j)} {f_2(n,i,j) \over g_2(n,i,j)},
\end{equation}
where $f_1$ and $g_1$ ($f_2$ and $g_2$, respectively) are
relatively prime polynomials.

Since the following discussion is independent of $n$, we omit the variable $n$
for convenience. For example, we write $R_1(i,j)$ instead of $R_1(n,i,j)$.
Using these notations, we have
\begin{theo}
\label{Zeil} Suppose the polynomials $g_1$ and $g_2$ in
\eqref{f-g} satisfy the hypotheses $(H1)$--$(H3)$. Suppose further
that for any $h_1,h_2 \in \{-1,0,1\}$,
\begin{equation}
\label{g-d} \gcd(g_1(i,j), d(i+h_1,j+h_2)) = \gcd(g_2(i,j),
d(i+h_1,j+h_2)) = 1.
\end{equation}

Then $g_1(i,j),g_2(i,j)$ can be factored into polynomials:
\[
\begin{array}{l}
g_1(i,j) = v_1(i) v_2(j) v_3(i+j) v_4(i,j) u_1(j) u_2(i,j), \\[5pt]
g_2(i,j) = v_1(i) v_2(j) v_3(i+j) v_4(i,j) w_1(i) w_2(i,j),
\end{array}
\]
such that
\[
\begin{array}{l}
v_1(i) \,|\, r_1(i-1,j) s_2'(i-1,j), \\[5pt]
v_2(j) \,|\, r_2(i,j-1) s_1'(i,j-1), \\[5pt]
v_4(i,j) \,|\, \gcd \big( r_1(i-1,j) s_2'(i-1,j), r_2(i,j-1) s_1'(i,j-1) \big), \\[5pt]
u_2(i,j) \,|\, \gcd \big( s_1(i,j) s_2'(i,j), r_1(i-1,j) s_2'(i-1,j) \big), \\[5pt]
w_2(i,j) \,|\, \gcd \big( s_2(i,j) s_1'(i,j), r_2(i,j-1) s_1'(i,j-1) \big),
\end{array}
\]
where
\begin{equation}
\label{s-s'}
\begin{array}{l}
s_1'(i,j)=s_1(i,j)/\gcd(s_1(i,j),s_2(i,j)), \\[5pt]
s_2'(i,j) = s_2(i,j)/ \gcd(s_1(i,j),s_2(i,j)).
\end{array}
\end{equation}

\end{theo}

\pf
Substituting \eqref{f-g} into \eqref{WZ-triple} and dividing $F(i,j)$ on both sides,
we obtain
\begin{eqnarray}
{c(i,j) \over d(i,j)} &=& {r_1(i,j) \over s_1(i,j)} {f_1(i+1,j) \over d(i+1,j)g_1(i+1,j)}
- {f_1(i,j) \over d(i,j)g_1(i,j)} \nonumber \\
&& + {r_2(i,j) \over s_2(i,j)} {f_2(i,j+1) \over d(i,j+1)g_2(i,j+1)}
- {f_2(i,j) \over d(i,j)g_2(i,j)} \label{xx},
\end{eqnarray}
i.e.,
\begin{eqnarray*}
c(i,j) &=& {r_1(i,j) d(i,j) \over s_1(i,j) d(i+1,j)} {f_1(i+1,j) \over g_1(i+1,j)}
- {f_1(i,j) \over g_1(i,j)} \\
&& + {r_2(i,j) d(i,j) \over s_2(i,j)d(i,j+1)} {f_2(i,j+1) \over g_2(i,j+1)}
- {f_2(i,j) \over g_2(i,j)}.
\end{eqnarray*}

Let
\[
\begin{array}{l}
\tilde{r}_1(i,j)=r_1(i,j) d(i,j), \quad \tilde{s}_1(i,j)=s_1(i,j) d(i+1,j), \\[5pt]
\tilde{r}_2(i,j)= r_2(i,j) d(i,j), \quad \tilde{s}_2(i,j)=s_2(i,j)d(i,j+1).
\end{array}
\]
All discussions in the proof of Theorem~\ref{Gosper} still hold.
Thus, we have
\begin{equation}
\label{tilde}
\begin{array}{l}
v_1(i) \,|\, \tilde{r}_1(i-1,j) \tilde{s}_2'(i-1,j), \\[5pt]
v_2(j) \,|\, \tilde{r}_2(i,j-1) \tilde{s}_1'(i,j-1), \\[5pt]
v_4(i,j) \,|\, \gcd \big( \tilde{r}_1(i-1,j) \tilde{s}_2'(i-1,j), \tilde{r}_2(i,j-1) \tilde{s}_1'(i,j-1) \big), \\[5pt]
u_2(i,j) \,|\, \gcd \big( \tilde{s}_1(i,j) \tilde{s}_2'(i,j), \tilde{r}_1(i-1,j) \tilde{s}_2'(i-1,j) \big), \\[5pt]
w_2(i,j) \,|\, \gcd \big( \tilde{s}_2(i,j) \tilde{s}_1'(i,j),
\tilde{r}_2(i,j-1) \tilde{s}_1'(i,j-1) \big),
\end{array}
\end{equation}
where
\[
\begin{array}{l}
\tilde{s}_1'(i,j) = \tilde{s}_1(i,j)/ \gcd( \tilde{s}_1(i,j),\tilde{s}_2(i,j) ), \\[5pt]
\tilde{s}_2'(i,j) = \tilde{s}_2(i,j)/ \gcd( \tilde{s}_1(i,j),\tilde{s}_2(i,j) ).
\end{array}
\]
Since we have \eqref{g-d}, we may replace
$\tilde{r}_1,\tilde{s}_1,\tilde{r}_2,\tilde{s}_2$ by $r_1,s_1,r_2,s_2$ in
\eqref{tilde}, respectively.  \qed

\section{\normalsize A Telescoping Method for Bivariate Hypergeometric Terms}

Theorem~\ref{Zeil} enables us to choose the denominators in the
telescoping algorithm. Basically, we will use certain  factors
appearing in the bounds of the denominators as estimates of the
denominators. In many cases, this approach seems to work quite
efficiently although we are not able to give a formula to  bound
the denominators because certain factors are not bounded in
Theorem~\ref{Zeil}. Roughly speaking, the divisibility
considerations in our method serve as a guide to guess the factors
in the denominators. In fact, the estimated denominators are much
smaller than the theoretical bounds given by Theorem~\ref{Zeil}.
Only $u_2(i,j)$ and $w_2(i,j)$ are set to their theoretical
bounds, while $v_2(j),v_3(i+j),v_4(i,j)$ are set to $1$, $u_1(j)$
and $w_1(i)$ are set to factors of $s_1(i,j)s'_2(i,j)$,  and
$v_1(i)$ is set to a factor of its theoretical bound. See the
following algorithm EstDen.

{\noindent \bf Algorithm EstDen} \\
{\bf Input:}  A hypergeometric term $F(n,i,j)$.\\
{\bf Output:} Estimated denominators $g_1(i,j)$ and $g_2(i,j)$ for
bivariate Gosper's algorithm.
\begin{itemize}
\item[1.] Calculate $r_1,r_2,s_1,s_2,s_1',s_2'$ defined by \eqref{IJ-rat} and
\eqref{s-s'}; \item[2.]
Set \\
$v_1(i):=$ the maximal factor of $r_1(i,j) s_2'(i,j)$ depending only on $i$; \\
$v_2(j):=$ the maximal factor of $r_2(i,j) s_1'(i,j)$ depending only on $j$; \\
and
\[
v(i):=\gcd(v_1(i-1),v_2(i-1));
\]
\item[3.]
Set \\
$u_1(j):=$ the maximal factor of $s_1(i,j)s'_2(i,j)$ depending only on $j$; \\
$w_1(i):=$ the maximal factor of $s_1(i,j)s'_2(i,j)$ depending only on $i$;
\item[4.] Set $u_2(i,j)$ to be the maximal factor of
\[
\gcd(s_1(i,j) s_2'(i,j), r_1(i-1,j)s_2'(i-1,j))
\]
which depends on $i$; \\
Set $w_2(i,j)$ to be the maximal factor of
\[
\gcd(s_1(i,j) s_2'(i,j), r_2(i,j-1)s_1'(i,j-1))
\]
which depends on $j$. \item[5.] Return $g_1(i,j):=v(i) u_1(j) u_2(i,j)$ and
$g_2(i,j):=v(i) w_1(i) w_2(i,j)$.
\end{itemize}

\noindent {\bf Remark.}
Let $f(i,j)$ be a polynomial in $i,j$ and $a$ be a new variable. Then the maximal
factor of $f(i,j)$ depending only on $i$ can be obtained by
\[
\gcd(f(i,j), f(i,j+a)),
\]
and the maximal factor of $f(i,j)$ depending on $i$
can be obtained by
\[
f(i,j)/\gcd(f(i,j), f(i+a,j)).
\]

We are now ready to describe our telescoping method for double
summations:

{\noindent \bf Method BiZeil}\\
{\bf Input:} A hypergeometric term $F(n,i,j)$.\\
{\bf Output:} An operator $L$ and rational functions $R_1$ and
$R_2$ such that \eqref{WZ-triple} holds if the algorithm succeeds.
\begin{itemize}
\item[1.] Using algorithm EstDen to obtain $g_1$ and $g_2$.
\item[2.] Set the order $r$ of the linear difference operator $L$ to be zero.
\item[3.] For the order $r$, calculate the common denominator $d(n,i,j)$ of
\[
{F(n+1,i,j) \over F(n,i,j)}, \quad \ldots, \quad {F(n+r,i,j) \over F(n,i,j)}.
\]
(If $r=0$, then take $d(n,i,j)=1$.)

\item[4.] Set the degrees of $f_1$ and $f_2$ to be one more than
those of $d \cdot g_1$ and $d \cdot g_2$, respectively.

\item[5.] Solve the equation \eqref{xx} by the method of
undeterminate coefficients to obtain $a_0,a_1,\ldots,a_r$ and
$f_1,f_2$.

\item[6.] If $a_i \not= 0$ for some $i \in \{0,\ldots,r\}$, then
return $L,f_1/(d \cdot g_1),f_2/(d \cdot g_2)$ and we are done.

If $a_i=0$ for all $i \in \{0,\ldots,r\}$, but $\deg f_1 - \deg (d
\cdot g_1) \le 2$, then increase the degrees of $f_1$ and $f_2$ by
one and repeat step 5.

Otherwise, set $r:=r+1$ and repeat the process from step 3.
\end{itemize}

\noindent {\bf Remarks.}
\begin{itemize}
\item[1.] In many cases, $g_1(i,j)$ and $g_2(i,j)$ can be further
reduced by cancelling a factor of degree $1$ and a factor of
degree $2$ from $g_1$ and $g_2$, respectively. In our
implementation we first choose two arbitrary factors and use the
reduced $g_1$ and $g_2$. When it fails, we then try the unreduced
ones. This cancellation may reduce the time of calculation if we
are lucky. For example, for the Andrews-Paule identity (see
Example $1$), the estimated denominators given by
Theorem~\ref{Zeil}, by algorithm EstDen, and by reduction are,
respectively, {\small
\begin{align*}
& g_1(i,j) = (2n-2i+1)(n-i+1)(2n-2j+1)(n-j+1)(i+j)^2 (j+1)^2, \\
& g_2(i,j) = (2n-2i+1)(n-i+1)(2n-2j+1)(n-j+1)(i+j)^2 (i+1)^2;
\\[5pt]
& g_1(i,j)=(2n-2i+1)(n-i+1)(j+1)^2, \\
& g_2(i,j)=(2n-2i+1)(n-i+1)(i+1)^2; \end{align*}} and
\[
g_1(i,j) = (2n-2i+1)(j+1)^2, \quad g_2(i,j)= (2n-2i+1)(n-i+1).
\]
The calculation times are $116$ seconds, $5$ seconds and $0.6$
second, respectively. We should note that since our method is
heuristic and it applies only to particular cases, we are more
interested in the computation results which are verifiable. So we
cannot claim the efficiency of the method or its applicability.

\item[2.] In all the following examples except Example $4$, the
degree of the numerator of $R_1$ ($R_2$) is one more than that of
the denominator. While in Example 4, the difference is two.

The degree bounds can be interpreted as follows. Let
$t_1,t_2,t_3,t_4$ be the four terms of the right hand side of
\eqref{xx} after multiplying the common denominator. In most
cases, the leading terms of $t_1$ and $t_2$ ($t_3$ and $t_4$,
respectively) are cancelled.

\item[3.] There is a way to speed up the computation in Step 5.
Given $g_1$ and $g_2$, we may derive part of the factors of $f_1$
and $f_2$ by divisibility. For example, suppose \eqref{xx} becomes
\[
{c(i,j) \over d(i,j)} = {u_1(i,j) \over v_1(i,j)} f_1(i+1,j)
- {f_1(i,j) \over w_1(i,j)} + {u_2(i,j) \over v_2(i,j)} f_2(i,j+1)
- {f_2(i,j) \over w_2(i,j)},
\]
after substituting and simplification. Suppose further that $D(i,j)$ is the common
denominator of the above equation. Then we immediately have that $f_1 \cdot
D/w_1$ is divisible by $q_1=\gcd(cD/d, u_1 D/v_1, u_2 D/v_2, D/w_2)$ and $f_1(i+1,j) \cdot
u_1D/v_1$ is divisible by $q_2=\gcd(cD/d, D/w_1, u_2D/v_2, D/w_2)$, and hence,
\[
{q_1 \over \gcd(D/w_1,q_1)} \quad \mbox{and} \quad
{q_2 \over \gcd(u_1D/v_1,q_2)}
\]
are factors of $f_1(i,j)$ and $f_1(i+1,j)$, respectively.
\end{itemize}

\section{\normalsize Examples}
In the following examples, let $F$ denote the summand of the left
hand side of the identity.

\noindent {\bf Example $1$.} The Andrews-Paule identity:
\begin{equation}
\label{A-P} \sum_{i=0}^n \sum_{j=0}^n {i+j \choose i}^2 {4n-2i-2j \choose
2n-2i} = (2n+1){2n \choose n}^2.
\end{equation}
It was proved by G.~Andrews and P.~Paule \cite{A-P92,A-P93} by
establishing a more general identity
\[
\sum_{i=0}^{\lfloor{m \over 2}\rfloor} \sum_{j=0}^{\lfloor{n \over 2}\rfloor}
{i+j \choose i}^2 {m+n-2i-2j \choose n-2i} = {\lfloor{m+n+1 \over 2}\rfloor!
\lfloor{m+n+2 \over 2}\rfloor! \over \lfloor{m \over 2}\rfloor!  \lfloor{m+1
\over 2}\rfloor!  \lfloor{n \over 2}\rfloor!  \lfloor{n+1 \over 2}\rfloor!}.
\]
Using the method BiZeil, we can deal with \eqref{A-P} directly. In
fact, we have
\[
g_1(i,j)=(2n-2i+1)(n-i+1)(j+1)^2, \quad g_2(i,j)=(2n-2i+1)(n-i+1)(i+1)^2.
\]
Cancelling the factors $(n-i+1)$ and $(i+1)^2$ from $g_1(i,j)$ and $g_2(i,j)$,
respectively, we obtain
\[
\tilde{g}_1(i,j) = (2n-2i+1)(j+1)^2
\quad \mbox{and} \quad \tilde{g}_2(i,j)= (2n-2i+1)(n-i+1).
\]

Finally, we get
\[
(2n+1) F(n,i,j) = \Delta_i R_1F(n,i,j) + \Delta_j R_2F(n,i,j),
\]
where
\[
R_1={i^2(6n^2+5n+1+6jn^2+jn-j-in+2in^2-2i-4j^2n-2j^2-3ij-4ijn) \over (2n-2i+1)(1+j)^2},
\]
\[
R_2 = {-2n^2+2jn^2+6in^2+9in+3jn-4ijn-4i^2n-n+j-3ij+2i-4i^2 \over (2n-2i+1)},
\]
which are the same as given in \cite[p.~85]{Weg}. Summing $i,j=0,\ldots,n$, we
get
\begin{eqnarray*}
\lefteqn{(2n+1) \sum_{i=0}^n \sum_{j=0}^n F(n,i,j)} \\
&=& \sum_{i=0}^n \big( R_2F(n,i,n+1)-R_2F(n,i,0) \big) + \sum_{j=0}^n \big(
R_1F(n,n+1,j)-R_1F(n,0,j) \big)
\end{eqnarray*}
Note that there is only one nonzero term $R_1F(n,n+1,n)$ of the
second summation. While applying Gosper's algorithm to the first
summand, we obtain
\[
\sum_{i=0}^n \big( R_2F(n,i,n+1)-R_2F(n,i,0) \big) = G(n+1)-G(0).
\]
where
\[
G(i) = {(-2n+i-1)(-4n+2i-1)i \over -1+2i-2n} {4n-2i \choose
2n-2i}.
\]
Simplifying $G(n+1)-G(0)+R_1F(n,n+1,n)$, we finally get
\eqref{A-P}.

\noindent {\bf Example $2$.} Carlitz's identity \cite{Car68} (see
Also \cite[Example 6.1.2]{W-Z92}):
\[
\sum_{i=0}^n \sum_{j=0}^n {i+j \choose i}{n-i \choose j}{n-j
\choose n-i-j} = \sum_{l=0}^n {2l \choose l}.
\]
We have
\[
g_1(i,j)=(j+1)^2(-n+j), \quad  g_2(i,j)=(i+1)^2(-n+i).
\]
Cancelling the factors $(-n+j)$ and $(i+1)(-n+i)$, we obtain
\[
\tilde{g}_1(i,j) = (j+1)^2
\quad \mbox{and} \quad \tilde{g}_2(i,j)= i+1.
\]
Notice that the common denominator of
\[
{F(n+1,i,j) \over F(n,i,j)} \quad \mbox{and} \quad {F(n+2,i,j)
\over F(n,i,j)}
\]
is $(-n+i-1+j)^2 (-n+i-2+j)^2$. We finally get
\[
L = (4n+6) -(8+5n)N + (n+2)N^2,
\]
and
\begin{multline*}
R_1=\big( -i^2(-n+i-1)(36-10ji^2n-13j^2ni+60j^2+60ji-2i^2-38j^2i\\
-8ji^2+10i^3+36n^3-11in^3-14jn^3-2i^4-92jn^2+8i^2n-80in+5j^2n^2\\
+8j^2i^2+88jin+42j^2n-172jn+24jin^2+5i^2n^2+3i^3n-54in^2+88n^2+4j^3n\\
-90j+6j^3-40i+5n^4+90n) \big) \Big/ \big( (-n+i-1+j)^2(-n+i-2+j)^2(j+1)^2 \big),
\end{multline*}
\vskip -30pt

\begin{multline*}
R_2=\big( (64-19ji^2n-6j^2ni+14j^2+74ji+54i^2-10j^2i-36ji^2+2i^3+39n^3\\
-16in^3-9jn^3-4i^4+6ji^3-53jn^2+50i^2n-176in+4j^2n^2+4j^2i^2+5n^4\\
+83jin+16j^2n-100jn+22jin^2+11i^2n^2+4i^3n-93in^2+112n^2-60j\\
-108i+140n)(-n-1+j) \big) \Big/ \big( (-n+i-2+j)^2(-n+i-1+j)^2
\big),
\end{multline*}
such that
\begin{equation}
\label{ex2} LF(n,i,j) = \Delta_i R_1F(n,i,j) + \Delta_j
R_2F(n,i,j),
\end{equation}

By summing \eqref{ex2} over $i,j$ from $0$ to $n$, one derives
that $L$ annihilates the double sum on the left hand side. It is
easily seen that the right hand side can be annihilated by $\big(
(n+2)N-(4n+6) \big)(N-1)$, which is exactly $L$. Then the identity
follows from  the initial values $n=0,1$.

The proofs of the following examples are similar to that of
Example $2$. We only need to give $\tilde{g}_1,\tilde{g}_2$ and
$L,R_1, R_2$. Then these identities can be verified by checking
the initial values.

\noindent {\bf Example $3$.} Carlitz's identity \cite{Car64} (see
also \cite[Example 6.1.3]{W-Z92}):
\begin{multline*}
\sum_{i=0}^m \sum_{j=0}^n {i+j \choose i}{m-i+j \choose j}{n-j+i
\choose i}{m+n-i-j \choose m-i}
 \\
= {(m+n+1)! \over m!n!} \sum_k {1 \over 2k+1} {m \choose k}{n \choose k}.
\end{multline*}
By cancelling the factors $(1+j)$ and $(i+1)^2$, we obtain
\[
\tilde{g}_1(i,j) = (n-j+i)(1+j) \quad \mbox{and} \quad
\tilde{g}_2(i,j)= m-i+j
\]
Notice that the common denominator of
\[
{F(n+1,i,j) \over F(n,i,j)} \quad \mbox{and} \quad {F(n+2,i,j)
\over F(n,i,j)}
\]
is $(-n+j-1)^2 (-n+j-2)^2$, which is denoted by $d(i,j)$. We
finally get
\begin{multline*}
L = 2(m+3+n)(2+m+n)^2 - \\
  (3m+2nm+4n^2+14+15n)(n+m+3)N +(2n+5)(n+2)^2 N^2,
\end{multline*}
and the denominators of $R_1,R_2$ are $d(i,j)\tilde{g}_1(i,j)$ and
$d(i,j)\tilde{g}_2(i,j)$, respectively. The degrees of
denominators and numerators of $R_1,R_2$ are both less than those
given in \cite{W-Z92}.

\noindent {\bf Example $4$.} The Ap\'ery-Schmidt-Strehl identity \cite{Stre}:
\[
\sum_i \sum_j {n \choose j} {n+j \choose j} {j \choose i}^3
= \sum_k {n \choose k}^2 {n+k \choose k}^2.
\]

By cancelling the factors $(-j-1+i)$ and $(i+1)^2$, we obtain
\[
\tilde{g}_1(i,j) = (-j-1+i)^2 \quad \mbox{and} \quad
\tilde{g}_2(i,j)= i+1
\]
Notice that the common denominator of
\[
{F(n+1,i,j) \over F(n,i,j)} \quad \mbox{and} \quad {F(n+2,i,j)
\over F(n,i,j)}
\]
is $(n+2-j)(n+1-j)$. We finally get
\[
L = (n+1)^3 - (3+2n)(17n^2+51n+39)N +(n+2)^3 N^2,
\]
and
\begin{multline*}
R_1  =  \big( -2i^2(3+2n)
 (-10+30j^2-49n^2-j^3-4n^4-24n^3-2n^2i^2+n^2i-6ni^2\\
 +3ni+3nji+n^2ji + 3j^2i^2-3j^3i+3ji-4i^2-2j^2i-2ji^2+11n^2j^2+6n^2j\\
 +33nj^2+18nj-6j^4+2i+15j-39n) \big) \Big/ \big( (n+2-j)(n+1-j)(-j-1+i)^2 \big),
\end{multline*}
\begin{multline*}
R_2  = \big(
2(-j+i)(3+2n)(-8n^2i-4n^2i^2-4n^2ji+4n^2j+4n^2j^2+12nj\\
-12nji-24ni+12nj^2-12ni^2+12j^2-4ji^2+j^3+6j^2i^2-3j^4+8j\\
+5j^2i-8i^2+3j^3i-16i-16ji) \big) \Big/ \big( (n+2-j)(n+1-j)(i+1) \big).
\end{multline*}
The rational functions $R_1,R_2$ are simpler than those given in \cite{Stre}.
The operator $L$ was used by Ap\'ery in his proof of the irrationality of
$\zeta(3)$ and Chyzak and Salvy obtained it using Ore algebras \cite{Chy-Sal}.

\noindent {\bf Example $5$.} The Strehl identity \cite{Stre}:
\begin{equation}
\label{complex}
\sum_i \sum_j {n \choose j} {n+j \choose j} {j \choose i}^2
{2i \choose i}^2 {2i \choose j-i} = \sum_k {n \choose k}^3 {n+k \choose k}^3.
\end{equation}
By cancelling the factor $(-3i-3+j)(-3i-2+j)$ from $g_2$, we
obtain
\[
\tilde{g}_1(i,j) = (j+1-i)^3
\quad \mbox{and} \quad \tilde{g}_2(i,j)= (-3i-1+j)(i+1)^3.
\]
Notice that the common denominator of
\[
{F(n+1,i,j) \over F(n,i,j)}, \ \ldots, \ {F(n+6,i,j) \over
F(n,i,j)}
\]
is $(n+1-j)(n+2-j)\cdots(n+6-j)$, which is denoted by $d(i,j)$. We
finally get a linear difference operator $L$ of order $6$ and the
denominators of $R_1,R_2$ are $d(i,j)\tilde{g}_1(i,j)$ and
$d(i,j)\tilde{g}_2(i,j)$, respectively. The operator $L$ is the
same as the operator obtained by applying Zeilberger's algorithm
to the right hand side of \eqref{complex}.

\noindent {\bf Example $6$.} The Graham-Knuth-Patashnik identity
\cite[p.~172]{Gra}:
\begin{equation}
\label{GKP} \sum_j \sum_k (-1)^{j+k}{j+k \choose k+l}{r \choose j}{n \choose
k}{s+n-j-k \choose m-j} = (-1)^l {n+r  \choose n+l}{s-r \choose m-n-l}.
\end{equation}
By cancelling the factor $(j+1)(j+1-l)$ from $g_2$, we obtain
\[
\tilde{g}_1(j,k) = (k+1)(k+l+1) \quad \mbox{and} \quad \tilde{g}_2(j,k)= 1.
\]
Notice that the denominator of ${F(r+1,j,k) \over F(r,j,k)}$ is
$r-j+1$, which is denoted by $d(j,k)$. We finally get a linear
difference operator with respect to the variable $r$:
\[
L=(r+n+1)(n+s+l-m-r)+(r-l+1)(r-s)R
\]
and the denominators of $R_1,R_2$ are $d(j,k)\tilde{g}_1(j,k)$ and
$d(j,k)\tilde{g}_2(j,k)$, respectively. Then \eqref{GKP} follows
from the evaluation of the initial value ($r=0$) by Zeilberger's
algorithm:
\[
\sum_k (-1)^k {k \choose k+l} {n \choose k}{s+n-k \choose m} = (-1)^l {n
\choose n+l}{s \choose m-n-l}.
\]

\noindent {\bf Example $7$.} The Petkov\v{s}ek-Wilf-Zeilberger identity
\cite[p.~33]{Pet-W-Z}:
\begin{equation}
\label{P-W-Z} \sum_r \sum_s (-1)^{n+r+s} {n \choose r} {n \choose s} {n+s
\choose s} {n+r \choose r} {2n-r-s \choose n} = \sum_k {n \choose k}^4.
\end{equation}
By cancelling the factors $s+1$ and $(r+1)^2$, we obtain
\[
\tilde{g}_1(r,s) = (n+r)(n+1-r)(s+1) \quad \mbox{and} \quad \tilde{g}_2(r,s)=
(n+r)(n+1-r).
\]
Notice that the common denominator of
\[
{F(n+1,r,s) \over F(n,r,s)} \quad \mbox{and} \quad {F(n+2,r,s)
\over F(n,r,s)}
\]
is
\[
(n+1)(n+2)(n+1-r)(n+2-r)(n+1-s)(n+2-s)(n-r-s+1)(n+2-r-s),
\]
which is denoted by $d(r, s)$. We finally get
\[
L=4(4n+5)(4n+3)(n+1)+2(2n+3)(3n^2+9n+7)N-(n+2)^3N^2
\]
is a linear difference operator and the denominators of $R_1,R_2$ are
$d(r,s)\tilde{g}_1(r,s)$ and $d(r,s)\tilde{g}_2(r,s)$, respectively. The
recursion is the same as that obtained by applying Zeilberger's algorithm to
the right hand side of \eqref{P-W-Z}.

\vskip 15pt {\small {\bf Acknowledgments.} We would like to thank
the referees for valuable suggestions leading to a significant
improvement of an earlier version. The authors are grateful to
Professor P. Paule for many helpful discussions and comments. This
work was supported by the ``973'' Project on Mathematical
Mechanization, the National Science Foundation, the Ministry of
Education, and the Ministry of Science and Technology of China.}

\end{document}